# q-Lucas polynomials and associated Rogers-Ramanujan type identities

Johann Cigler
Fakultät für Mathematik, Universität Wien
johann.cigler@univie.ac.at

**Abstract**
We prove some properties of $q-$ analogues of the Fibonacci and Lucas polynomials, apply these to derive some identities due to L. Carlitz und H. Prodinger and finally give an easy approach to L. Slater's Bailey pairs A(1)-A(8) using $q-$Lucas polynomials.

## 1. Introduction

After recalling some properties of the $q-$Fibonacci and $q-$Lucas polynomials which I have introduced in [4] and [5] I apply these to derive some identities due to L. Carlitz [2] und H. Prodinger [8]. Finally I show that the $q-$Lucas polynomials allow an easy approach to L. Slater's Bailey pairs A(1)-A(8) and some related Rogers-Ramanujan type identities. I want to thank Andrew Sills for commenting on a previous version and pointing out to me the papers [10] and [11] of L. Slater.

Let

$$F_n(x,s) = \sum_{k=0}^{\frac{n-1}{2}} \binom{n-1-k}{k} s^k x^{n-1-2k} \qquad (1.1)$$

and

$$L_n(x,s) = \sum_{k=0}^{\frac{n}{2}} \frac{n}{n-k}\binom{n-k}{k} s^k x^{n-2k} \qquad (1.2)$$

be the classical Fibonacci and Lucas polynomials. They satisfy the recurrence

$$F_n(x,s) = xF_{n-1}(x,s) + sF_{n-2}(x,s) \qquad (1.3)$$

with initial values $F_0(x,s) = 0, F_1(x,s) = 1$ and

$$L_n(x,s) = xL_{n-1}(x,s) + sL_{n-2}(x,s) \qquad (1.4)$$

with initial values $L_0(x,s) = 2$ and $L_1(x,s) = x$.

It will be convenient to define a variant $L_n^*(x,s)$ by $L_n^*(x,s) = L_n(x,s)$ for $n > 0$ and $L_0^*(x,s) = 1$.

Let $\alpha = \dfrac{x+\sqrt{x^2+4s}}{2}$ and $\beta = \dfrac{x-\sqrt{x^2+4s}}{2}$ be the roots of the equation $z^2 - xz - s = 0$.

Then it is well-known and easily verified that $F_n(x,s) = \dfrac{\alpha^n - \beta^n}{\alpha - \beta}$ and $L_n(x,s) = \alpha^n + \beta^n$.
This implies the well-known formula

$$L_n(x,s) = F_{n+1}(x,s) + sF_{n-1}(x,s) \tag{1.5}$$

for $n > 0$ because $s = \alpha\beta$ and $\alpha^{n+1} - \beta^{n+1} - \alpha\beta(\alpha^{n-1} - \beta^{n-1}) = \alpha^n(\alpha - \beta) + \beta^n(\alpha - \beta)$.

Another known formula is

$$\sum_{k=0}^{\lfloor \frac{n}{2} \rfloor} (-s)^k \binom{n}{k} L^*_{n-2k}(x,s) = x^n. \tag{1.6}$$

For the proof it is convenient to consider this identity for odd $n$ and even $n$ separately.
For odd $n$ the left-hand side is

$$\sum_{k=0}^{\lfloor \frac{n}{2} \rfloor}(\alpha\beta)^k \binom{n}{k}\left(\alpha^{n-2k} + \beta^{n-2k}\right) = \sum_{k=0}^{\lfloor \frac{n}{2} \rfloor}\binom{n}{k}\left(\alpha^{n-k}\beta^k + \alpha^k\beta^{n-k}\right) = \sum_{k=0}^{n}\binom{n}{k}\alpha^k\beta^{n-k} = (\alpha+\beta)^n = x^n.$$

For $n = 2m$ the same holds because for $k = m$ the coefficient of $(-s)^m \binom{2m}{m}$ is $L^*_0(x,s) = 1$.

Using (1.5) we see that (1.6) is equivalent with

$$\sum_{k=0}^{\lfloor \frac{n}{2} \rfloor}\left(\binom{n}{k} - \binom{n}{k-1}\right)(-s)^k F_{n+1-2k}(x,s) = x^n. \tag{1.7}$$

Let us first review the simpler case of $q$-analogues of the binomial theorem
$(x+s)^n = \sum_{k=0}^{n}\binom{n}{k} x^{n-k} s^k = (x+s)(x+s)^{n-1}$. It is well known that there are two important ones,

$$p_n(x,s) = (x+s)(x+qs)\cdots(x+q^{n-1}s) = \sum_{k=0}^{n} q^{\binom{k}{2}} \begin{bmatrix} n \\ k \end{bmatrix} x^{n-k} s^k, \tag{1.8}$$

which satisfy the recurrence relation $p_n(x,s) = (x+q^{n-1}s)p_{n-1}(x,s)$

and the Rogers-Szegö polynomials

$$r_n(x,s) = \sum_{k=0}^{n} \begin{bmatrix} n \\ k \end{bmatrix} x^{n-k} s^k, \tag{1.9}$$

which have no closed formula but satisfy the recursion

$$r_n(x,s) = (x+s)r_{n-1}(x,s) + \left(q^{n-1}-1\right)sxr_{n-2}(x,s). \tag{1.10}$$

Here $\begin{bmatrix} n \\ k \end{bmatrix} = \dfrac{(q;q)_n}{(q;q)_k(q;q)_{n-k}}$ with $(x;q)_n = \prod_{j=0}^{n-1}\left(1-q^j x\right)$ denotes the $q$-binomial coefficient. We also use $[n]$ instead of $\begin{bmatrix} n \\ 1 \end{bmatrix}$.

A similar situation occurs with $q$-analogues of the Fibonacci polynomials. There are the polynomials studied by L. Carlitz

$$f_n(x,s,q) = \sum_{k=0}^{\lfloor \frac{n-1}{2} \rfloor} \begin{bmatrix} n-1-k \\ k \end{bmatrix} q^{k^2} x^{n-1-2k} s^k, \tag{1.11}$$

which satisfy the recursion $f_n(x,s,q) = xf_{n-1}(x,qs,q) + qsf_{n-2}(x,q^2s,q)$
and the polynomials $F_n(x,s,q)$ with which I am concerned in this paper.

## 2. Definition and simple properties

Define the $q$-Fibonacci polynomials $F_n(x,s,q)$ by

$$F_n(x,s,q) = \sum_{k=0}^{\frac{n-1}{2}} q^{\binom{k+1}{2}} \begin{bmatrix} n-1-k \\ k \end{bmatrix} s^k x^{n-1-2k} \tag{2.1}$$

for $n \geq 0$.

The first polynomials are

$0, 1, x, x^2 + qs, x^3 + (1+q)qsx, x^4 + qs[3]x^2 + q^3s^2, \cdots$

Let us recall that these $q$-Fibonacci polynomials satisfy each of the recurrences

$$F_n(x,s,q) = xF_{n-1}(x,qs,q) + qsF_{n-2}(x,qs,q), \tag{2.2}$$

$$F_n(x,s,q) = xF_{n-1}(x,s,q) + q^{n-2}sF_{n-2}(x,\frac{s}{q},q) \tag{2.3}$$

and

$$F_n(x,s,q) = xF_{n-1}(x,s,q) + q^{n-2}sxF_{n-3}(x,s,q) + q^{n-2}s^2F_{n-4}(x,s,q). \tag{2.4}$$

The simple proofs follow by comparing coefficients and using the well-known recurrences for the $q-$binomial coefficients. We see that (2.2) is equivalent with

$$\begin{bmatrix} n-1-k \\ k \end{bmatrix} = q^k \begin{bmatrix} n-2-k \\ k \end{bmatrix} + \begin{bmatrix} n-2-k \\ k-1 \end{bmatrix}$$

and (2.3) with

$$\begin{bmatrix} n-1-k \\ k \end{bmatrix} = \begin{bmatrix} n-2-k \\ k \end{bmatrix} + q^{n-1-2k} \begin{bmatrix} n-2-k \\ k-1 \end{bmatrix}.$$

Combining (2.2) and (2.3) we get (2.4).

As a consequence we get
$$F_n(x,s,q) - xF_{n-1}(x,s,q) - sF_{n-2}(x,s,q)$$
$$= -sF_{n-2}(x,s,q) + q^{n-2}sF_{n-2}(x,s,q) - q^{n-2}sF_{n-2}(x,s,q) + q^{n-2}sxF_{n-3}(x,s,q) + q^{n-2}s^2F_{n-4}(x,s,q)$$
$$= (q^{n-2}-1)sF_{n-2}(x,s,q) - q^{n-2}s(F_{n-2}(x,s,q) - xF_{n-3}(x,s,q) - sF_{n-4}(x,s,q))$$

Iterating this equation and observing that it holds for $n=2$ and $n=3$ gives

$$F_n(x,s,q) - xF_{n-1}(x,s,q) - sF_{n-2}(x,s,q) = \sum_{k=1}^{\lfloor n/2 \rfloor} (-1)^k q^{(k-1)(n-k)} s^k \left(1 - q^{n-2k}\right) F_{n-2k}(x,s,q). \qquad (2.5)$$

**Remark 1**
From (2.2) we get the following combinatorial interpretation of the $q-$Fibonacci polynomials which is a $q-$analogue of the well-known Morse code model of the Fibonacci numbers. Consider words $c = c_1 c_2 \cdots c_m$ of letters $c_i \in \{a,b\}$ and associate with $c$ the weight $w(c) = w(c)(s) = q^{i_1 + i_2 + \cdots + i_k} s^k x^{m-k}$, if $c_{i_1} = c_{i_2} = \cdots = c_{i_k} = b, 1 \le i_1 < \cdots < i_k \le m$, and all other $c_i = a$. The weight of the empty word $\varepsilon$ is defined to be $w(\varepsilon) = 1$.
We then have
$$\begin{aligned} w(ac)(s) &= xw(c)(qs), \\ w(bc)(s) &= qsw(c)(qs), \\ w(ca)(s) &= xw(c)(s), \\ w(cb)(s) &= q^{m+1}sw(c)(s). \end{aligned} \qquad (2.6)$$

Define the length $l(c)$ of a word $c$ consisting of $k$ letters $b$ and $m-k$ letters $a$ by $l(c) = 2k + m - k = m + k$. Let now $\Phi_n$ be the set of all words of length $n-1$. $\Phi_n$ can also be identified with the set of all coverings of an $(n-1) \times 1-$rectangle with monominos (i.e. $1 \times 1-$rectangles) and dominos (i.e. $2 \times 1-$rectangles) or with Morse code sequences of length $n-1$.

Let $G_n(x,s) := \sum_{c \in \Phi_n} w(c)$ be the weight of $\Phi_n$.

Then we get
$$G_n(x,s) := \sum_{c \in \Phi_n} w(c) = F_n(x,s,q), \tag{2.7}$$

i.e. $F_n(x,s,q)$ is the weight of $\Phi_n$.

For the proof observe that by considering the first letter of each word we see from (2.6) that $G_n(x,s) = xG_{n-1}(x,qs) + qsG_{n-2}(x,qs)$. Thus $G_n(x,s)$ satisfies the same recurrence as $F_n(x,s,q)$. Also the initial values coincide because $G_0(x,s) = 0$ and $G_1(x,s) = 1$.

For example $\Phi_4 = \{aaa, ab, ba\}$ and
$$G_4(x,s) = w(aaa) + w(ab) + w(ba) = x^3 + xq^2s + qsx = F_4(x,s,q).$$

This interpretation can also be used to obtain (2.3), which is equivalent with
$$G_n(x,s) = xG_{n-1}(x,s) + q^{n-2}sG_{n-2}\left(x, \frac{s}{q}\right).$$

Here we consider the last letter of each word. From (2.6) we get $w(ca)(s) = xw(c)(s)$, which gives the first term. To obtain the second term let us suppose that $cb \in \Phi_n$ has $k$ letters $b$. Then $w(cb)(s) = q^m sw(c)(s) = q^{n-1-k}sq^{i_1 + \cdots + i_{k-1}}s^{k-1}x^{n-k-1} = q^{n-2}sw(c)\left(\frac{s}{q}\right)$. Since this expression is independent of $k$, we get the second term.

Let $D$ be the $q$-differentiation operator defined by $Df(x) = \dfrac{f(x) - f(qx)}{(1-q)x}$. As has been shown in [5] these $q$-Fibonacci polynomials also satisfy
$$F_n(x,s,q) = F_n\big(x + (q-1)sD, s\big)1. \tag{2.8}$$

In order to show (2.8) we must verify that
$$F_n(x,s,q) = xF_{n-1}(x,s,q) + (q-1)sDF_{n-1}(x,s,q) + sF_{n-2}(x,s,q). \tag{2.9}$$

Comparing coefficients this amounts to
$$q^k \begin{bmatrix} n-1-k \\ k \end{bmatrix} = q^k \begin{bmatrix} n-2-k \\ k \end{bmatrix} + \begin{bmatrix} n-2-k \\ k-1 \end{bmatrix} + \left(q^{n-2k} - 1\right)\begin{bmatrix} n-1-k \\ k-1 \end{bmatrix}$$
or
$$\left(q^k - 1\right)\begin{bmatrix} n-1-k \\ k \end{bmatrix} = \left(q^{n-2k} - 1\right)\begin{bmatrix} n-1-k \\ k-1 \end{bmatrix}$$
which is obviously true.

As in [5] we define the $q$-Lucas polynomials by

$$L_n(x,s,q) = L_n\big(x+(q-1)sD, s\big)1. \tag{2.10}$$

The first polynomials are

$$2, x, x^2+(1+q)s, x^3+[3]sx, x^4+[4]sx^2+q(1+q^2)s^2, \cdots$$

By applying the linear map

$$f(x) \to f(x+(q-1)sD)1 \tag{2.11}$$

to (1.5) we get

$$L_n(x,s,q) = F_{n+1}(x,s,q) + sF_{n-1}(x,s,q) \tag{2.12}$$

for $n > 0$.

This implies the explicit formula

$$L_n(x,s,q) = \sum_{k=0}^{\frac{n}{2}} q^{\binom{k}{2}} \frac{[n]}{[n-k]} \begin{bmatrix} n-k \\ k \end{bmatrix} s^k x^{n-2k} \tag{2.13}$$

for $n > 0$, which is a very nice $q$-analogue of (1.2).

For the proof observe that

$$q^k \begin{bmatrix} n-k \\ k \end{bmatrix} + \begin{bmatrix} n-k-1 \\ k-1 \end{bmatrix} = \frac{q^k[n-k]+[k]}{[n-k]} \begin{bmatrix} n-k \\ k \end{bmatrix} = \frac{[n]}{[n-k]} \begin{bmatrix} n-k \\ k \end{bmatrix}.$$

Comparing coefficients we also get

$$L_n(x,qs,q) = F_{n+1}(x,s,q) + q^n s F_{n-1}(x,s,q). \tag{2.14}$$

This follows from

$$q^k \begin{bmatrix} n-k \\ k \end{bmatrix} + q^n \begin{bmatrix} n-k-1 \\ k-1 \end{bmatrix} = \frac{q^k[n-k]+q^n[k]}{[n-k]} \begin{bmatrix} n-k \\ k \end{bmatrix} = q^k \frac{[n-k]+q^{n-k}[k]}{[n-k]} \begin{bmatrix} n-k \\ k \end{bmatrix}$$

$$= q^k \frac{[n]}{[n-k]} \begin{bmatrix} n-k \\ k \end{bmatrix}.$$

**Remark 2**
(2.12) has the following combinatorial interpretation:
Consider a circle whose circumference has length $n$ and let monominos be arcs of length 1 and dominos be arcs of length 2 on the circle. Consider the set $\Lambda_n$ of all coverings with monominos and dominos and fix a point $P$ on the circumference of the circle. If $P$ is the initial point of a monomino or a domino of a covering then this covering can be identified with a word $c = c_1 \cdots c_m$. We define its weight in the same way as in the linear case. Therefore the set of all those coverings has weight $F_{n+1}(x,s,q)$.

If $P$ is the midpoint of a domino we split $b$ into $b = b_0 b_1$ and associate with this covering the word $b_1 c_1 \cdots c_m b_0$ with $c_1 \cdots c_m \in \Phi_{n-1}$ and define its weight as $sw(c_1 \cdots c_m)$.

Therefore $w(\Lambda_n) = w(\Phi_{n+1}) + sw(\Phi_{n-1})$.

E.g. $\Lambda_4 = \{aaaa, aab, aba, baa, bb, b_1 aab_0, b_1 bb_0\}$. Thus

$$w(\Lambda_4) = x^4 + x^2 q^3 s + xq^2 sx + qsx^2 + qsq^2 s + sx^2 + sqs = x^4 + [4]sx^2 + q(1+q^2)s^2 = L_4(x,s,q).$$

To give a combinatorial interpretation of (2.14) we consider all words of $\Lambda_n$ with last letter $a$ or the two last letters $ab_0$. Their weight is $F_n(x,s,q) + sF_{n-1}(x,s,q) = F_{n+1}\left(x, \dfrac{s}{q}, q\right)$.

There remains the set of all words in $\Lambda_n$ with last letter $b$. With the same argument as above we see that this is $q^{n-1} s F_{n-1}\left(x, \dfrac{s}{q}, q\right)$. Therefore we have

$$L_n(x,s,q) = F_{n+1}\left(x, \dfrac{s}{q}, q\right) + q^{n-1} s F_{n-1}\left(x, \dfrac{s}{q}, q\right),$$

which is equivalent with (2.14).

We also need the polynomials $L_n^*(x,s,q)$ which coincide with $L_n(x,s,q)$ for $n > 0$, but have initial value $L_0^*(x,s,q) = 1$.

Comparing (2.10) with (1.4) we see that

$$L_n(x,s,q) = xL_{n-1}(x,s,q) + (q-1)sDL_{n-1}(x,s,q) + sL_{n-2}(x,s,q). \tag{2.15}$$

This is a recurrence for the polynomials in $x$ but not for individual numbers $x$ and $s$.

In order to find a recurrence for individual numbers I want to show first that for $n > 2$

$$L_n^*(x,s,q) - xL_{n-1}^*(x,s,q) - sL_{n-2}^*(x,s,q) = (1-q^{n-1})\sum_{k=1}^{\lfloor n/2 \rfloor}(-1)^k q^{(k-1)(n-k-1)} s^k L_{n-2k}^*(x,s,q). \tag{2.16}$$

This reduces to (1.4) for $q = 1$.

It is easily verified that $DL_n(x,s,q) = [n]F_n\left(x, \dfrac{s}{q}, q\right)$.

Therefore

$$L_n(x,s,q) - xL_{n-1}(x,s,q) - sL_{n-2}(x,s,q) = (q-1)sDL_{n-1}(x,s,q) = \left(q^{n-1}-1\right)sF_{n-1}\left(x, \dfrac{s}{q}, q\right).$$

By (2.14) we know that $F_{n-1}\left(x,\frac{s}{q},q\right) = L_{n-2}(x,s,q) - q^{n-3}sF_{n-3}\left(x,\frac{s}{q},q\right).$
Iteration gives (2.16).

From (2.16) we get
$L_n(x,s,q) - xL_{n-1}(x,s,q) - sL_{n-2}(x,s,q) = (q^{n-1}-1)sL_{n-2}(x,s,q)$
$-q^{n-3}s\frac{[n-1]}{[n-3]}(L_{n-2}(x,s,q) - xL_{n-3}(x,s,q) - sL_{n-4}(x,s,q))$

This can be written as
$L_n(x,s,q) = xL_{n-1}(x,s,q) - \frac{(1+q)q^{n-3}}{[n-3]}sL_{n-2}(x,s,q) + \frac{[n-1]}{[n-3]}q^{n-3}sxL_{n-3}(x,s,q) + \frac{[n-1]}{[n-3]}q^{n-3}s^2L_{n-4}(x,s,q)$ (2.17)

This recurrence holds for $n \geq 4$ if $L_0(x,s,q) = 2$.

### 3. Inversion theorems

L. Carlitz [2] has obtained two $q-$ analogues of the Chebyshev inversion formulas. The first one ([2],Theorem [6]) implies

**Theorem 3.1**

$$\sum_{2k \leq n} \begin{bmatrix} n \\ k \end{bmatrix} L^*_{n-2k}(x,s,q)(-s)^k = x^n \qquad (3.1)$$

and the second one ([2], Theorem 7) gives

**Theorem 3.2**

$$\sum_{2k \leq n} \left(\begin{bmatrix} n \\ k \end{bmatrix} - \begin{bmatrix} n \\ k-1 \end{bmatrix}\right) F_{n+1-2k}(x,s,q)(-s)^k = x^n. \qquad (3.2)$$

These are $q-$ analogues of (1.6) and (1.7). We give another proof of these theorems:

Let
$$A = x + (q-1)sD. \qquad (3.3)$$

Define
$$\alpha(q) = \frac{A + \sqrt{A^2 + 4s}}{2} \qquad (3.4)$$

and

$$\beta(q) = \frac{A - \sqrt{A^2 + 4s}}{2}. \qquad (3.5)$$

Since $\alpha(q)^2 - A\alpha(q) - s = \beta(q)^2 - A\beta(q) - s = 0$ the sequences $\left(\alpha(q)^n\right)_{-\infty}^{\infty}$ and $\left(\beta(q)^n\right)_{-\infty}^{\infty}$ satisfy the recurrence

$$\alpha(q)^n - A\alpha(q)^{n-1} - s\alpha(q)^{n-2} = \beta(q)^n - A\beta(q)^{n-1} - s\beta(q)^{n-2} = 0$$

for all $n \in \mathbb{Z}.$

Since the $q-$Fibonacci and the $q-$Lucas polynomials satisfy the same recurrence we get from the initial values

$$L_n(x,s,q) = \left(\alpha(q)^n + \beta(q)^n\right)1 \qquad (3.6)$$

and

$$F_n(x,s,q) = \frac{\alpha(q)^n - \beta(q)^n}{\alpha(q) - \beta(q)} 1 \qquad (3.7)$$

for $n \geq 0.$ We can use these identities to extend these polynomials to negative $n$.

We then get for $n > 0$

$$L_{-n}(x,s,q) = \left(\alpha(q)^{-n} + \beta(q)^{-n}\right)1 = (-1)^n \frac{\beta(q)^n + \alpha(q)^n}{s^n} = (-1)^n \frac{L_n(x,s,q)}{s^n} \qquad (3.8)$$

and

$$F_{-n}(x,s,q) = \frac{\alpha(q)^{-n} - \beta(q)^{-n}}{\alpha(q) - \beta(q)} 1 = (-1)^{n-1} \frac{F_n(x,s,q)}{s^n}. \qquad (3.9)$$

**Remark 3**
It is easily verified that the identities (2.2), (2.3), (2.4), (2.12) and (2.14) hold for all $n \in \mathbb{Z}.$

We want to show that

$$\sum_{k=0}^{\left\lfloor \frac{n}{2} \right\rfloor} (-s)^k \begin{bmatrix} n \\ k \end{bmatrix} L^*_{n-2k}(x,s,q) = x^n. \qquad (3.10)$$

For odd $n$ the left-hand side is

$$\sum_{k=0}^{\lfloor \frac{n}{2} \rfloor}(\alpha(q)\beta(q))^k \begin{bmatrix} n \\ k \end{bmatrix}\left(\alpha(q)^{n-2k}+\beta(q)^{n-2k}\right)1 = \sum_{k=0}^{\lfloor \frac{n}{2} \rfloor}\begin{bmatrix} n \\ k \end{bmatrix}\left(\alpha(q)^{n-k}\beta(q)^k + \alpha(q)^k\beta(q)^{n-k}\right)1 = \sum_{k=0}^{n}\begin{bmatrix} n \\ k \end{bmatrix}\alpha(q)^k\beta(q)^{n-k}1.$$

For $n=2m$ the same holds because for $k=m$ the coefficient of $(-s)^m \begin{bmatrix} 2m \\ m \end{bmatrix}$ is $L_0^*(x,s,q)=1$.

Therefore by (1.10) we see that $R(n,x,s) := \sum_{k=0}^{\lfloor \frac{n}{2} \rfloor}(-s)^k \begin{bmatrix} n \\ k \end{bmatrix} L_{n-2k}^*(x,s,q)$ satisfies

$$R(n,x,s) = (\alpha(q)+\beta(q))R(n-1,x,s) + \alpha(q)\beta(q)\left(q^{n-1}-1\right)R(n-2,x,s)$$
$$= (x+(q-1)sD)R(n-1,x,s) - s\left(q^{n-1}-1\right)R(n-2,x,s).$$

We have to show that $R(n,x,s) = x^n$.
This is obviously true for $n=0$ and $n=1$.
If it holds for $m < n$ then

$$R(n,x,s) = (x+(q-1)sD)x^{n-1} - s\left(q^{n-1}-1\right)x^{n-2} = x^n$$

as asserted.

### Remark 4

In [6] we have defined a new $q-$analogue of the Hermite polynomials
$H_n(x,s\,|\,q) = (x-sD)^n 1$. By applying the linear map (2.11) to (1.6) and (1.7) these can be expressed as

$$H_n(x,(q-1)s\,|\,q) = \sum_{k=0}^{n/2}\binom{n}{k}s^k L_{n-2k}^*(x,-s,q) = \sum_{k=0}^{(n+1)/2}\left(\binom{n}{k} - \binom{n}{k-1}\right)s^k F_{n+1-2k}(x,-s,q). \quad (3.11)$$

### Remark 5

The polynomials $\left(F_{n+1}(x,s,q)\right)_{n \geq 0}$ are a basis of the polynomials in $\mathbb{C}(s,q)[x]$. Define a linear functional $L$ on this vector space by $L(F_{n+1}(x,s,q)) = [n=0]$.
Then from (3.2) we get

$L(x^{2n+1}) = 0$ and $L(x^{2n}) = (-qs)^n C_n(q)$, where $C_n(q) = \dfrac{1}{[n+1]}\begin{bmatrix} 2n \\ n \end{bmatrix}$ is a $q-$analogue of the Catalan numbers. This is equivalent with

$$\sum_{k=0}^{n}(-1)^k q^{\binom{k}{2}}\begin{bmatrix} 2n-k \\ k \end{bmatrix}C_{n-k}(q) = [n=0]. \quad (3.12)$$

In the same way the linear functional $M$ defined by $M(L_n^*(x,s,q)) = [n=0]$ gives $M(x^{2n}) = \begin{bmatrix} 2n \\ n \end{bmatrix}(-s)^n$ and $M(x^{2n+1}) = 0$. This is equivalent with

$$\sum_{k=0}^{n}(-1)^k q^{\binom{k}{2}} \frac{[2n]}{[2n-k]} \begin{bmatrix} 2n-k \\ k \end{bmatrix} \begin{bmatrix} 2n-2k \\ n-k \end{bmatrix} = 0 \tag{3.13}$$

for $n > 0$.

## 4. Some related identities

The classical Fibonacci and Lucas polynomials satisfy

$$L_n(x+y, -xy) = x^n + y^n \tag{4.1}$$

and

$$F_n(x+y, -xy) = \frac{x^n - y^n}{x - y}. \tag{4.2}$$

L. Carlitz [2] has given $q-$analogues of these theorems which are intimately connected with our $q-$analogues.

**Theorem 4.1 (L. Carlitz[2])**

Let $r_n(x,y) = \sum_{k=0}^{n} \begin{bmatrix} n \\ k \end{bmatrix} x^k y^{n-k}$ be the Rogers-Szegö polynomials. Then

$$\frac{x^n - y^n}{x - y} = \sum_{k=0}^{\lfloor \frac{n-1}{2} \rfloor} q^{\binom{k+1}{2}} \begin{bmatrix} n-k-1 \\ k \end{bmatrix} (-xy)^k r_{n-1-2k}(x,y) \tag{4.3}$$

and

$$x^n + y^n = \sum_{k=0}^{\lfloor \frac{n}{2} \rfloor} q^{\binom{k}{2}} \frac{[n]}{[n-k]} \begin{bmatrix} n-k \\ k \end{bmatrix} (-xy)^k r_{n-2k}(x,y). \tag{4.4}$$

These are polynomial identities which for $(x,y) \to (\alpha, \beta)$ immediately give the explicit formulae for the $q-$Fibonacci and $q-$Lucas-polynomials if we define them by (3.7) and (3.6).

To prove these theorems we use the identity

$$\sum_{j=0}^{k}(-1)^j q^{\binom{j+1}{2}} \begin{bmatrix} k \\ j \end{bmatrix} \begin{bmatrix} n-j \\ k \end{bmatrix} = 1 \tag{4.5}$$

for $0 \leq k \leq n$ (cf. Carlitz[2]).

To show this identity let $n$ be fixed and let $U$ be the $\mathbb{C}(q)$-linear operator on the vector space of all sums $\sum_{k=0}^{n} c_k \begin{bmatrix} x-k \\ n-k \end{bmatrix}$ with $c_k \in \mathbb{C}(q)$ defined by $U \begin{bmatrix} x-k \\ n-k \end{bmatrix} = \begin{bmatrix} x-k-1 \\ n-k \end{bmatrix}$ for $0 \leq k \leq n$.

Since $(1 - q^{n-k} U) \begin{bmatrix} x-k \\ n-k \end{bmatrix} = \begin{bmatrix} x-k \\ n-k \end{bmatrix} - q^{n-k} \begin{bmatrix} x-k-1 \\ n-k \end{bmatrix} = \begin{bmatrix} x-k-1 \\ n-k-1 \end{bmatrix}$

by using (1.8) we get the desired result

$$\sum_{j=0}^{n} (-1)^j q^{\binom{j+1}{2}} \begin{bmatrix} n \\ j \end{bmatrix} \begin{bmatrix} x-j \\ n \end{bmatrix} = (1-qU)\cdots(1-q^n U) \begin{bmatrix} x \\ n \end{bmatrix} = (1-qU)\cdots(1-q^{n-1}U) \begin{bmatrix} x-1 \\ n-1 \end{bmatrix} = \cdots = \begin{bmatrix} x-n \\ 0 \end{bmatrix}.$$

First we prove (4.3).

$$\sum_{k=0}^{\lfloor \frac{n-1}{2} \rfloor} q^{\binom{k+1}{2}} \begin{bmatrix} n-k-1 \\ k \end{bmatrix} (-xy)^k r_{n-1-2k}(x,y) = \sum_{k=0}^{\lfloor \frac{n-1}{2} \rfloor} q^{\binom{k+1}{2}} \begin{bmatrix} n-k-1 \\ k \end{bmatrix} (-xy)^k \sum_{j=0}^{n-1-2k} \begin{bmatrix} n-1-2k \\ j \end{bmatrix} x^j y^{n-1-2k-j}$$

$$= \sum_{i=0}^{n-1} x^i y^{n-1-i} \sum_{k=0}^{i} (-1)^k q^{\binom{k+1}{2}} \begin{bmatrix} n-k-1 \\ k \end{bmatrix} \begin{bmatrix} n-1-2k \\ i-k \end{bmatrix} = \sum_{i=0}^{n-1} x^i y^{n-1-i} \sum_{k=0}^{i} (-1)^k q^{\binom{k+1}{2}} \begin{bmatrix} i \\ k \end{bmatrix} \begin{bmatrix} n-1-k \\ i \end{bmatrix}$$

$$= \sum_{i=0}^{n-1} x^i y^{n-1-i} \begin{bmatrix} n-1-i \\ 0 \end{bmatrix} = \sum_{i=0}^{n-1} x^i y^{n-1-i} = \frac{x^n - y^n}{x-y}.$$

For the proof of (4.4) observe that

$$\frac{x^{n+1} - y^{n+1}}{x-y} - xy \frac{x^{n-1} - y^{n-1}}{x-y} = \frac{x^n(x-y) + y^n(x-y)}{x-y} = x^n + y^n.$$

This implies

$$x^n + y^n = \sum_k q^{\binom{k+1}{2}} \begin{bmatrix} n-k \\ k \end{bmatrix} (-xy)^k r_{n-2k}(x,y) - xy \sum_k q^{\binom{k+1}{2}} \begin{bmatrix} n-k-2 \\ k \end{bmatrix} (-xy)^k r_{n-2-2k}(x,y)$$

$$= \sum_k q^{\binom{k+1}{2}} \begin{bmatrix} n-k \\ k \end{bmatrix} (-xy)^k r_{n-2k}(x,y) - \sum_k q^{\binom{k}{2}} \begin{bmatrix} n-k-1 \\ k-1 \end{bmatrix} (-xy)^k r_{n-2k}(x,y)$$

$$= \sum_k q^{\binom{k}{2}} \left( q^k \begin{bmatrix} n-k \\ k \end{bmatrix} - \begin{bmatrix} n-k-1 \\ k-1 \end{bmatrix} \right) (-xy)^k r_{n-2k}(x,y) = \sum_{k=0}^{\lfloor \frac{n}{2} \rfloor} q^{\binom{k}{2}} \frac{[n]}{[n-k]} \begin{bmatrix} n-k \\ k \end{bmatrix} (-xy)^k r_{n-2k}(x,y).$$

For the classical Fibonacci polynomials the formula

$$\sum_{k=0}^{n}(-1)^k x^k \binom{n}{k} F_{2n+m-k}(x,s) = s^n F_m(x,s) \tag{4.6}$$

holds for all $m \in \mathbb{Z}$. This is an easy consequence of the Binet formula $F_n(x,s) = \dfrac{\alpha^n - \beta^n}{\alpha - \beta}$,

where $\alpha = \dfrac{x + \sqrt{x^2 + 4s}}{2}$ and $\beta = \dfrac{x - \sqrt{x^2 + 4s}}{2}$. For (4.6) is equivalent with

$$\frac{\alpha^{n+m}(\alpha - x)^n - \beta^{n+m}(\beta - x)^n}{\alpha - \beta} = s^n \frac{\alpha^m - \beta^m}{\alpha - \beta}.$$

We now get

**Theorem 4.2**

$$\sum_{k=0}^{n}(-1)^k q^{\binom{k}{2}} \begin{bmatrix} n \\ k \end{bmatrix} x^k F_{2n+m-k}(x,s,q) = q^{\binom{n}{2}+mn} s^n F_m\left(x, \frac{s}{q^n}, q\right). \tag{4.7}$$

The case $m = 0$ gives

**Corollary 4.3 (H. Prodinger [8])**

$$\sum_{k=0}^{n}(-1)^k q^{\binom{k}{2}} \begin{bmatrix} n \\ k \end{bmatrix} x^k F_{2n-k}(x,s,q) = 0. \tag{4.8}$$

**Proof**
For $n = 0$ this is trivially true for all $m \in \mathbb{Z}$.
For $n = 1$ (4.7) reduces to

$$F_{m+2}(x,s,q) - x F_{m+1}(x,s,q) = q^m s F_m\left(x, \frac{s}{q}, q\right), \tag{4.9}$$

which also holds for $m \in \mathbb{Z}$ by (2.3) and Remark 3.

Assume that (4.7) holds for $i < n$ and all $m$. Then we get

$$\sum_k (-1)^k q^{\binom{k}{2}} \begin{bmatrix} n \\ k \end{bmatrix} x^k F_{2n+m-k}(x,s,q) = \sum_k (-1)^k q^{\binom{k}{2}} q^k \begin{bmatrix} n-1 \\ k \end{bmatrix} x^k F_{2n+m-k}(x,s,q)$$

$$+ \sum_k (-1)^k q^{\binom{k}{2}} \begin{bmatrix} n-1 \\ k-1 \end{bmatrix} x^k F_{2n+m-k}(x,s,q)$$

$$= \sum_k (-1)^{k-1} q^{\binom{k}{2}} \begin{bmatrix} n-1 \\ k-1 \end{bmatrix} x^{k-1} F_{2n+m-k+1}(x,s,q) + \sum_k (-1)^k q^{\binom{k}{2}} \begin{bmatrix} n-1 \\ k-1 \end{bmatrix} x^k F_{2n+m-k}(x,s,q)$$

$$= \sum_k (-1)^{k-1} q^{\binom{k}{2}} \begin{bmatrix} n-1 \\ k-1 \end{bmatrix} x^{k-1} \left( F_{2n+m-k+1}(x,s,q) - x F_{2n+m-k}(x,s,q) \right)$$

$$= \sum_k (-1)^k q^{\binom{k+1}{2}} \begin{bmatrix} n-1 \\ k \end{bmatrix} x^k q^{2n+m-k-2} s F_{2n+m-k-2}\left(x, \frac{s}{q}, q\right)$$

$$= q^{2n+m-2} s \sum_k (-1)^k q^{\binom{k}{2}} \begin{bmatrix} n-1 \\ k \end{bmatrix} x^k F_{2(n-1)+m-k}\left(x, \frac{s}{q}, q\right)$$

$$= q^{2n+m-2} q^{\binom{n-1}{2}+m(n-1)} s \left(\frac{s}{q}\right)^{n-1} F_m\left(x, \frac{s}{q^n}, q\right) = q^{\binom{n}{2}+mn} s^n F_m\left(x, \frac{s}{q^n}, q\right).$$

For

$$h(n,m) = \sum_{k=0}^{n} (-1)^k q^{\binom{k}{2}} \begin{bmatrix} n \\ k \end{bmatrix} x^k L_{2n+m-1-k}(x,s,q) \tag{4.10}$$

we get from

$$L_n(x,s,q) = F_{n+1}(x,s,q) + s F_{n-1}(x,s,q) \tag{4.11}$$

$$h(n,m) = q^{\binom{n}{2}+mn} s^n F_m\left(x, \frac{s}{q^n}, q\right) + s^{n+1} q^{\binom{n}{2}+mn-2n} F_{m-2}\left(x, \frac{s}{q^n}, q\right). \tag{4.12}$$

This implies

$$\sum_{k=0}^{n} (-1)^k q^{\binom{k}{2}} \left( q^{n-1} \begin{bmatrix} n \\ k \end{bmatrix} + \begin{bmatrix} n-1 \\ k \end{bmatrix} \right) L_{2n+m-1-k}(x,s,q) = q^{n-1} h(n,m) + h(n-1, m+2).$$

Because of (3.9) we get

$$\sum_{k=0}^{n} (-1)^k q^{\binom{k}{2}} \left( q^{n-1} \begin{bmatrix} n \\ k \end{bmatrix} + \begin{bmatrix} n-1 \\ k \end{bmatrix} \right) L_{2n-1-k}(x,s,q) = q^{n-1} h(n,0) + h(n-1, 2)$$

$$= q^{n-1} \left( s^{n+1} q^{\binom{n}{2}-2n} F_{-2}\left(x, \frac{s}{q^n}, q\right) \right) + q^{\binom{n-1}{2}+2(n-1)} s^{n-1} F_2\left(x, \frac{s}{q^{n-1}}, q\right)$$

$$= -s^{n+1} q^{\binom{n}{2}-n-1} \left(\frac{q^n}{s}\right)^2 x + q^{\binom{n-1}{2}+2n-2} s^{n-1} x = 0.$$

This gives

**Theorem 4.4 (H. Prodinger [8] )**

$$\sum_{k=0}^{n}(-1)^k q^{\binom{k}{2}}\left(q^{n-1}\begin{bmatrix}n\\k\end{bmatrix}+\begin{bmatrix}n-1\\k\end{bmatrix}\right)x^k L_{2n-1-k}(x,s,q)=0. \quad (4.13)$$

## 5. Some Rogers-Ramanujan type formulas

It is interesting that the $q$ – Lucas polynomials give a simple approach to the Bailey pairs A(1) - A(8) of Slater's paper [10].

Let us recall some definitions (cf. [1] or [7] ) suitably modified for our purposes.
Two sequences $a=(\alpha_n)$ and $b=(\beta_n)$ are called a Bailey pair $(a,b)_m$, if

$$\beta_n = \sum_{k=0}^{n}\frac{\alpha_k}{(q;q)_{n-k}(q;q)_{n+k+m}} \quad (5.1)$$

for some $m \in \{0,1\}$. Note that $(q;q)_n = (1-q)(1-q^2)\cdots(1-q^n)$.

To obtain Bailey pairs we start with formula (3.1) and consider separately even and odd numbers $n$. This gives

$$\sum_{k=0}^{n}\begin{bmatrix}2n\\n-k\end{bmatrix}L^*_{2k}(x,-s(q),q)s(q)^{n-k} = x^{2n} \quad (5.2)$$

and

$$\sum_{k=0}^{n}\begin{bmatrix}2n+1\\n-k\end{bmatrix}L^*_{2k+1}(x,-s(q),q)s(q)^{n-k} = x^{2n+1}. \quad (5.3)$$

Therefore

**Theorem 5.1**

$$a = \left(L^*_{2n}(x,-s(q),q)s(q)^{-n}\right), b = \left(\frac{x^{2n}}{s(q)^n (q;q)_{2n}}\right) \quad (5.4)$$

*is a Bailey pair with* $m=0$

*and*

$$a = \left(L^*_{2n+1}(x,-s(q),q)s(q)^{-n}\right), b = \left(\frac{x^{2n+1}}{s(q)^n (q;q)_{2n+1}}\right) \quad (5.5)$$

*one with* $m=1$.

If we change $q \to \frac{1}{q}$ we get the Bailey pairs

$$a = \left(L_{2n}^*(x, -s(q^{-1}), q^{-1})q^{n^2} s(q^{-1})^{-n}\right), b = \left(\frac{x^{2n} q^{n^2}}{s(q^{-1})^n (q;q)_{2n}}\right) \quad (5.6)$$

with $m = 0$ and

$$a = \left(L_{2n+1}^*(x, -s(q^{-1}), q^{-1})q^{n^2+n} s(q^{-1})^{-n}\right), b = \left(\frac{x^{2n+1} q^{n^2+n}}{s(q^{-1})^n (q;q)_{2n+1}}\right) \quad (5.7)$$

with $m = 1$.

For each Bailey pair we consider the identity

$$\sum_{n \geq 0} q^{n^2+mn} \beta_n = \sum_{n \geq 0} q^{n^2+mn} \sum_{k=0}^{n} \frac{\alpha_k}{(q;q)_{n-k}(q;q)_{n+k+m}} = \sum_{k \geq 0} \alpha_k \sum_{n \geq k} \frac{q^{n^2+mn}}{(q;q)_{n-k}(q;q)_{n+k+m}}. \quad (5.8)$$

Here the inner sum $\sum_{n \geq k} \frac{q^{n^2+nm}}{(q;q)_{n-k}(q;q)_{n+k+m}}$ can be easily computed:

For $k \in \mathbb{N}$ we have

$$\sum_{s \geq i} \frac{q^{s^2+ks}}{(q;q)_{s-i}(q;q)_{s+i+k}} = \frac{q^{i^2+ki}}{(q;q)_\infty} \quad (5.9)$$

This is an easy consequence of the $q-$Vandermonde formula

$$\sum_{s=i}^{n-i} q^{(s-i)(s+i+k)} \begin{bmatrix} n+2i \\ s+i+k \end{bmatrix} \begin{bmatrix} n-2i \\ s-i \end{bmatrix} = \sum_{j=0}^{n-2i} q^{j(2i+j+k)} \begin{bmatrix} n+2i \\ n-j-k \end{bmatrix} \begin{bmatrix} n-2i \\ j \end{bmatrix} = \begin{bmatrix} 2n \\ n-k \end{bmatrix}$$

if we let $n \to \infty$.

Therefore we get

$$\sum_{n \geq 0} q^{n^2+mn} \beta_n = \sum_{k \geq 0} \alpha_k \sum_{n \geq k} \frac{q^{n^2+mn}}{(q;q)_{n-k}(q;q)_{n+k+m}} = \frac{1}{(q;q)_\infty} \sum_{k \geq 0} \alpha_k q^{k^2+mk}. \quad (5.10)$$

In the following formulas we set $x = 1$ and $m \in \{0,1\}$.

For $s = 1$ we get from (5.10), (5.4) and (5.5)

$$\sum_s \frac{q^{n^2+mn}}{(q;q)_{2n+m}} = \frac{1}{(q;q)_\infty} \sum_{i \geq 0} L_{2i+m}^*(1, -1, q) q^{i^2+mi}. \quad (5.11)$$

For $s = \dfrac{1}{q}$ we get

$$\sum_{n \geq 0} \frac{q^{n^2+n+mn}}{(q;q)_{2n+m}} = \frac{1}{(q;q)_\infty} \sum_{k \geq 0} L^*_{2k+m}(1,-q^{-1},q) q^{k^2+k+mk}. \tag{5.12}$$

In the same way we get from (5.6) and (5.7) for $s = 1$

$$\sum_{s \geq 0} \frac{q^{2n^2+2mn}}{(q;q)_{2n+m}} = \sum_{s \geq 0} q^{n^2+mn} \sum_{i=0}^{s} \frac{L^*_{2i+m}(1,-1,q^{-1}) q^{i^2}}{(q;q)_{n-i}(q;q)_{n+i+m}} = \frac{1}{(q;q)_\infty} \sum_{i \geq 0} L^*_{2i+m}(1,-1,q^{-1}) q^{2i^2+2mi} \tag{5.13}$$

and for $s(q) = \dfrac{1}{q}$

$$\sum_{n \geq 0} \frac{q^{2n^2-n+2mn}}{(q;q)_{2n+m}} = \sum_{n \geq 0} q^{n^2+mn} \sum_{i=0}^{n} \frac{L^*_{2i+m}(1,-q,q^{-1}) q^{i^2-i+mi}}{(q;q)_{n-i}(q;q)_{n+i+m}} = \frac{1}{(q;q)_\infty} \sum_{i \geq 0} L^*_{2i+m}(1,-q,q^{-1}) q^{2i^2-i+2mi}. \tag{5.14}$$

The main advantage of these formulas derives from the fact, that the $q$-Lucas polynomials have simple values for $x = 1$ and $s = -1$ or $s = -\dfrac{1}{q}$.

From (2.4) it is easily verified (cf. [5]) that

$$F_{3n}\left(1,-\frac{1}{q},q\right) = 0, \quad F_{3n+1}\left(1,-\frac{1}{q},q\right) = (-1)^n q^{\frac{n(3n-1)}{2}}, \quad F_{3n+2}\left(1,-\frac{1}{q},q\right) = (-1)^n q^{\frac{n(3n+1)}{2}}. \tag{5.15}$$

Therefore by (2.14)

$$L_{3n}(1,-1,q) = (-1)^n \left( q^{\frac{n(3n-1)}{2}} + q^{\frac{n(3n+1)}{2}} \right) \text{ for n>0}$$

$$L_{3n+1}(1,-1,q) = (-1)^n q^{\frac{n(3n+1)}{2}}, \tag{5.16}$$

$$L_{3n-1}(1,-1,q) = (-1)^n q^{\frac{n(3n-1)}{2}}.$$

Of course in all formulas $L^*_0(1,s,q) = 1$, although I shall not state this in each case explicitly. (5.16) implies

$$L^*_{6n-2}(1,-1,q) = L^*_{2(3n-1)}(1,-1,q) = -q^{6n^2-5n+1}, \quad L^*_{6n}(1,-1,q) = L^*_{2(3n)}(1,-1,q)$$
$$= q^{6n^2-n} + q^{6n^2+n}, \quad L^*_{6n+2}(1,-1,q) = L^*_{2(3n+1)}(1,-1,q) = -q^{6n^2+5n+1} \tag{5.17}$$

and

$$L^*_{6n-1}(1,-1,q) = L^*_{2(3n-1)+1}(1,-1,q) = q^{6n^2-n}, \quad L^*_{6n+3}(1,-1,q) = L^*_{2(3n+1)+1}(1,-1,q)$$
$$= -\left( q^{6n^2+5n+1} + q^{6n^2+7+2} \right), \quad L^*_{6n+1}(1,-1,q) = L^*_{2(3n+1)-1}(1,-1,q) = q^{6n^2+n}. \tag{5.18}$$

The first terms of the sequence $L^*_{2n}(1,-1,q)$ are therefore
$1, -q, -q^2, q^5+q^7, -q^{12}, -q^{15}, q^{22}+q^{26}, \cdots$.
The sum of all these terms is Euler's pentagonal number series.
The same is true for the sequence $L^*_{2n+1}(1,-1,q)$, which begins with
$1, -q-q^2, q^5, q^7, -q^{12}-q^{15}, q^{22}, q^{26}, \cdots$.
This is an immediate consequence of (5.4) and (5.5) for $x = s = 1$, which reduce to

$$\sum_{k=0}^{n} \frac{L^*_{2k}(1,-1,q)}{(q;q)_{n-k}(q;q)_{n+k}} = \frac{1}{(q;q)_{2n}} \tag{5.19}$$

and

$$\sum_{k=0}^{n} \frac{L^*_{2k+1}(1,-1,q)}{(q;q)_{n-k}(q;q)_{n+k+1}} = \frac{1}{(q;q)_{2n+1}}. \tag{5.20}$$

If we let $n \to \infty$ these formulas converge to $\sum_{k=0}^{\infty} L^*_{2k}(1,-1,q) = (q;q)_\infty$
and $\sum_{k=0}^{\infty} L^*_{2k+1}(1,-1,q) = (q;q)_\infty$ respectively.

By (2.12) we get

$$L_{3n}(1,-q^{-1},q) = (-1)^n \left( q^{\frac{n(3n-1)}{2}} + q^{\frac{n(3n-5)}{2}} \right) \text{ for n>0}$$

$$L_{3n+1}(1,-q^{-1},q) = (-1)^n q^{\frac{n(3n+1)}{2}}, \tag{5.21}$$

$$L_{3n-1}(1,-q^{-1},q) = (-1)^n q^{\frac{(n-2)(3n-1)}{2}}.$$

This implies that

$$L^*_{6n-2}\left(1,-\frac{1}{q},q\right) = L^*_{2(3n-1)}\left(1,-\frac{1}{q},q\right) = -q^{6n^2-5n+1}, L^*_{6n}\left(1,-\frac{1}{q},q\right) = L^*_{2(3n)}\left(1,-\frac{1}{q},q\right)$$
$$= q^{6n^2-n} + q^{6n^2-5n}, L^*_{6n+2}\left(1,-\frac{1}{q},q\right) = L^*_{2(3n+1)}\left(1,-\frac{1}{q},q\right) = -q^{6n^2-n-1} \tag{5.22}$$

and

$$L^*_{6n-1}\left(1,-\frac{1}{q},q\right) = L^*_{2(3n-1)+1}\left(1,-\frac{1}{q},q\right) = q^{6n^2-7n+1}, L^*_{6n+3}\left(1,-\frac{1}{q},q\right) = L^*_{2(3n+1)+1}\left(1,-\frac{1}{q},q\right)$$
$$= -q^{6n^2+5n+1} - q^{6n^2+n-1}, L^*_{6n+1}\left(1,-\frac{1}{q},q\right) = L^*_{2(3n)+1}\left(1,-\frac{1}{q},q\right) = q^{6n^2+n}. \tag{5.23}$$

Now it is time to harvest the Corollaries. We order them so that Corollary 5.i corresponds to Slater's $A(i)$.

**Corollary 5.1** (cf. [9], A.79)

$$\sum_{n\geq 0}\frac{q^{n^2}}{(q;q)_{2n}}=\frac{1}{(q;q)_\infty}\sum_{k\in\mathbb{Z}}\left(q^{15k^2+k}-q^{15k^2+11k+2}\right). \tag{5.24}$$

**Proof.**
Choose $m=0$ in (5.11) and observe that
$L_{6i}(1,-1,q)q^{(3i)^2}=q^{15i^2-i}+q^{15i^2+i}$, $L_{6i+2}(1,-1,q)q^{(3i+1)^2}=q^{15i^2+11i+2}$ and
$L_{6i-2}(1,-1,q)q^{(3i-1)^2}=q^{15i^2-11i+2}$
which implies
$\sum_{i\geq 0}L^*_{2i}(1,-1,q)q^{i^2}=\sum_{k\in\mathbb{Z}}\left(q^{15k^2+k}-q^{15k^2+11k+2}\right)$ and thus (5.24).

**Corollary 5.2** (cf. [9], A.94)

$$\sum_{n\geq 0}\frac{q^{n^2+n}}{(q;q)_{2n+1}}=\frac{1}{(q;q)_\infty}\sum_{k\in\mathbb{Z}}\left(q^{15k^2-4k}-q^{15k^2+14k+3}\right). \tag{5.25}$$

**Proof.**
We use formula (5.11) for $m=1$ and compute

$L^*_{6n-1}(1,-1,q)q^{(3n-1)^2+3n-1}=q^{15n^2-4n}$, $L^*_{6n+3}\left(1,-\frac{1}{q},q\right)q^{(3n+1)^2+3n+1}=-q^{15n^2+14n+3}-q^{15n^2+16n+4}$,

$L^*_{6n+1}\left(1,-\frac{1}{q},q\right)q^{(3n)^2+3n}=q^{15n^2+4n}$.

Since $15(n-1)^2+16(n-1)+4=15n^2-14n+3$ we get (5.25).

**Corollary 5.3** (cf. [9], A.99)

$$\sum_{n\geq 0}\frac{q^{n^2+n}}{(q;q)_{2n}}=\frac{1}{(q;q)_\infty}\sum_{k\geq 0}L^*_{2k}(1,-q^{-1},q)q^{k^2+k}=\frac{1}{(q;q)_\infty}\sum_{k\in\mathbb{Z}}\left(q^{15k^2+2k}-q^{15k^2+8k+1}\right). \tag{5.26}$$

This follows from (5.12) for $m = 0$ and

$$L^*_{6n-2}\left(1, -\frac{1}{q}, q\right) q^{(3n-1)^2+3n-1} = -q^{15n^2-8n+1}, \quad L^*_{6n}\left(1, -\frac{1}{q}, q\right) q^{(3n)^2+3n} = q^{15n^2+2n} + q^{15n^2-2n},$$

$$L^*_{6n+2}\left(1, -\frac{1}{q}, q\right) q^{(3n+1)^2+3n+1} = -q^{15n^2+8n+1}.$$

**Corollary 5.4** (cf. [9], A.38)

$$\sum_{n \geq 0} \frac{q^{n^2+2n}}{(q;q)_{2n+1}} = \frac{1}{(q;q)_\infty} \sum_{k \in \mathbb{Z}} \left( q^{15k^2-7k} - q^{15k^2+13k+2} \right). \tag{5.27}$$

**Proof.**

This follows from (5.12) for $m = 1$ and the computation

$$L^*_{6n-1}\left(1, -\frac{1}{q}, q\right) q^{(3n-1)^2+2(3n-1)} = q^{15n^2-7n}, \quad L^*_{6n+3}\left(1, -\frac{1}{q}, q\right) q^{(3n+1)^2+2(3n+1)} = -q^{15n^2+17n+4} - q^{15n^2+13n+2},$$

$$L^*_{6n+1}\left(1, -\frac{1}{q}, q\right) q^{(3n)^2+2(3n)} = q^{15n^2+7n}.$$

Observe that $15(n-1)^2 + 17(n-1) + 4 = 15n^2 - 13n + 2$.

**Corollary 5.5** (cf. [9], A.39)

$$\sum_{n \geq 0} \frac{q^{2n^2}}{(q;q)_{2n}} = \frac{1}{(q;q)_\infty} \sum_{k \in \mathbb{Z}} \left( q^{12k^2+k} - q^{12k^2+7k+1} \right). \tag{5.28}$$

**Proof.**
Here we use (5.13) with $m = 0$.

$$L^*_{6n-2}\left(1, -1, \frac{1}{q}\right) q^{2(3n-1)^2} = -q^{12n^2-7n+1}, \quad L^*_{6n}\left(1, -1, \frac{1}{q}\right) q^{2(3n)^2} = q^{12n^2+n} + q^{12n^2-n},$$

$$L^*_{6n+2}\left(1, -1, \frac{1}{q}\right) q^{2(3n+1)^2} = -q^{12n^2+7n+1}.$$

**Corollary 5.6** (cf. [9], A.84)

$$\sum_{n\geq 0}\frac{q^{2n^2+n}}{(q;q)_{2n+1}}=\frac{1}{(q;q)_\infty}\sum_{k\in\mathbb{Z}}(-1)^k q^{3k^2-k}=\frac{(q^2;q^2)_\infty}{(q;q)_\infty}=(-q;q)_\infty. \quad (5.29)$$

**Proof.**
We use (5.14) with $m=1$.

$$L^*_{6n-1}\left(1,-q,\frac{1}{q}\right)q^{2(3n-1)^2+3n-1}=q^{3(2n)^2-(2n)}, L^*_{6n+3}\left(1,-q,\frac{1}{q}\right)q^{2(3n+1)^2+3n+1}=-q^{3(2n+1)^2-(2n+1)}-q^{3(2n+1)^2+(2n+1)},$$

$$L^*_{6n+1}\left(1,-q,\frac{1}{q}\right)q^{2(3n)^2+3n}=q^{3(2n)^2+(2n)}.$$

Therefore we get
$$\sum_i L^*_{2i+1}(1,-q,q^{-1})q^{2i^2+i}=\sum_{k\in\mathbb{Z}}(-1)^k q^{3k^2-k}.$$

**Corollary 5.7** (cf. [9], A.52)

$$\sum_{n\geq 0}\frac{q^{2n^2-n}}{(q;q)_{2n}}=\frac{1}{(q;q)_\infty}\sum_{k\in\mathbb{Z}}(-1)^k q^{3k^2-k}=\frac{(q^2;q^2)_\infty}{(q;q)_\infty}=(-q;q)_\infty. \quad (5.30)$$

**Proof.**
This follows from (5.14) with $m=0$, because we get the same sums as in Corollary 5.6.

$$L^*_{6n-2}\left(1,-q,\frac{1}{q}\right)q^{2(3n-1)^2+(3n-1)}=-q^{3(2n+1)^2-(2n+1)}, L^*_{6n}\left(1,-q,\frac{1}{q}\right)q^{2(3n)^2-(3n)}=q^{3(2n)^2-2n}+q^{3(2n)^2+2n},$$

$$L^*_{6n+2}\left(1,-q,\frac{1}{q}\right)q^{2(3n+1)^2-(3n+1)}=-q^{3(2n+1)^2-(2n+1)}.$$

The deeper reason for the simple results (5.29) and (5.30) are the formulas

$$L^*_{2i+m}\left(1,-q,\frac{1}{q}\right)q^{2i^2-i+2mi}=L^*_{2i+m}\left(1,-1,q^2\right) \quad (5.31)$$

for $m\in\{0,1\}$, which can easily be verified.

**Corollary 5.8** (cf. [9], A.96)

$$\sum_{n\geq 0}\frac{q^{2n^2+2n}}{(q;q)_{2n+1}} = \frac{1}{(q;q)_\infty}\sum_{k\in\mathbb{Z}}\left(q^{12k^2+5k} - q^{12k^2-13k+3}\right). \tag{5.32}$$

**Proof.**
Here we use (5.13) with $m=1$.
We get
$$L^*_{6n-1}\left(1,-1,\frac{1}{q}\right)q^{2(3n-1)^2+2(3n-1)} = q^{12n^2-5n}, \quad L^*_{6n+3}\left(1,-1,\frac{1}{q}\right)q^{2(3n+1)^2+2(3n+1)} = -q^{12n^2+13n+3} - q^{12n^2+11n+2},$$
$$L^*_{6n+1}\left(1,-1,\frac{1}{q}\right)q^{2(3n)^2+2(3n)} = q^{12n^2+5n}.$$

We have only to verify that $12(n-1)^2 + 11(n-1) + 2 = 12n^2 - 13n + 3$.